\documentclass[12pt,a4paper]{article}
\usepackage{amsfonts}
\usepackage{mathrsfs}
\usepackage{amsmath}
\usepackage{amssymb}
\usepackage{booktabs}

\renewcommand{\baselinestretch}{1}
\renewcommand{\thesection}

\begin{document}
\renewcommand\refname{\centerline{References}}
\title{{On $\big(p_{1}(x), p_{2}(x)\big)$-Laplace Equations
\thanks{Research supported by the National Natural Science Foundation of China
(NSFC 10971088 and NSFC 10971087).}}}
\author{Duchao Liu\footnote{Email address: liuduchao@yahoo.cn}, Xiaoyan Wang\footnote{Email address: wang264@indiana.edu} and Jinghua Yao\footnote{Email address: yaoj@indiana.edu}
 \\
{\small School of Mathematics and Statistics, Lanzhou
University, Lanzhou, 730000, China} \\
{\small Department of Mathematics, Indiana University Bloomington, IN, 47405, USA} \\
\small Department of Mathematics, Indiana University Bloomington,
IN, 47405, USA}
\date{}
\maketitle

\begin{center}
\begin{minipage}{120mm}
\baselineskip 12pt {\small {\normalsize\textbf{Abstract.}} In this
paper, we investigate the  $(p_{1}(x), p_{2}(x))$-Laplace operator,
the properties of the corresponding integral functional and weak
solutions to the related differential equations. We show that the
integral functional admits a derivative of type $(S_+)$ which induces
a homeomorphism between duality space pairs. As applications of the
above results, we gave some existence results of the $(p_{1}(x),
p_{2}(x))$-Laplace equation
\begin{align*}
-\text{div}(|\nabla u|^{p_{1}(x)-2}\nabla u)-\text{div}(|\nabla u|^{p_{2}(x)-2}\nabla u)=f(x,u)
\end{align*}
in a bounded smooth domain $\Omega\subset\mathbb{R}^{N}$ with
Dirichlet boundary condition.\\[6pt]
}
\end{minipage}
\end{center}

\section*{0. Introduction}

In the mathematical modeling of electrorheological fluids (see \cite{1,2,3} 
and the references therein), the $p(x)$-Laplace operator,
defined as $\Delta_{p(x)}u(x):=\text{div}\big(|\nabla
u(x)|^{p(x)-2}\nabla u\big)$, plays an important role. Based on
these application backgrounds, the study of differential equations
involving $p(x)$-Laplace operators have been a very interesting and
exciting topic in recent years (see in particular the nice survey
\cite{4,5} and the references therein for further details). During the
development of last several decades, it turns out that the use of
variational methods in dealing with problems involving nonstandard
growth conditions is a far-reaching field. The $p(x)-$growth
conditions can be regarded as an important particular case of the
nonstandard $(p, q)-$growth conditions. Many results have been
obtained on this kind of problems, for example
\cite{6,7,8,9,10,11,12,13} and so on.

In this paper, continuing our former investigations on these topics
(see \cite{6,7,8,11,12,13,14,15}), we shall study the joint effects
of different $(p_{1}(x),p_{2}(x))$-Laplace operators and investigate the properties
of the corresponding integral functional in proper framework of
variable exponent Lebesgue-Sobolev spaces (see Section 1). We choose
to present here these properties as it is clear that not only these
results have their own interests from the point view of pure
functional analysis but also they have direct applications in the
study of differential equation with $(p_1 (x), p_2 (x))$-growth
conditions.

As an illustration of the aforementioned properties, we choose to
present some results on the existence of weak solutions to the
following $(p_{1}(x), p_{2}(x))$-Laplace equation
\begin{align*}
(P)\left\{ \begin{array}{rcl}
         &-\text{div}(|\nabla u|^{p_{1}(x)-2}\nabla u)-\text{div}(|\nabla u|^{p_{2}(x)-2}\nabla u)=f(x,u), &  \mbox{ in }
         \Omega; \\
         &u=0,&\mbox{ on } \partial\Omega,
          \end{array}\right.
\end{align*}
where $\Omega\subset\mathbb{R}^{N}$ is a bounded smooth domain and
$p_{i}(x)\in C(\overline{\Omega})$ with $p_{i}(x)>1$ for any
$x\in\overline{\Omega}$ and for $i=1, 2$. These results themselves
are interesting and new. Under proper growth conditions on nonlinear
terms, we showed the corresponding functional enjoys coercive
property (Theorem 3.3), mountain pass geometry (Theorem 3.6) and
fountain geometry (Theorem 3.7) respectively, which yields rich
existence results to our problem. The role of the results on $(p_1
(x), p_2 (x))$-Laplace operator lies in that they supply us a neat
means to prove the required compactness result (Lemma 3.5) when we
use variational devices. It is important to notice that the growth
conditions on the nonlinear term are given according to the function
$p_M (x)$ and $p_m (x)$ other than $p_1^- \wedge p_2^-$ or $p_1^+
\vee p_2^+$ (see Section 3).

For the ease of exposition and to keep the paper in a reasonable
length, we only give existence results on equations with $(p_1 (x),
p_2 (x))$-growth conditions on \textit{bounded domain} and with
\textit{Dirichlet boundary condition}. With some symmetry conditions
(see \cite{16}) and weighted exponent Sobolev spaces method, we could also
give many results on unbounded domain with other types of boundary
conditions. Besides, under proper conditions, the functionals will
satisfy some other geometry structures. We will address these issues
in a following paper.

This paper is organized as follows. In Section 1, for the
convenience of the readers, we recall the definitions of the variable
exponent Lebesgue-Sobolev spaces which can be regarded as a special
class of generalized Orlicz-Sobolev spaces and introduce some basic
properties of these spaces; In Section 2, we investigate properties
of the $(p_{1}(x), p_{2}(x))$-Laplace operator and the corresponding
integral functional; In Section 3, we show the existence of weak
solutions to problem $(P)$ by variational arguments.

\section*{1. The spaces $W_{0}^{1,p(x)}(\Omega)$}
In this section, we will give out some theories on spaces
$W_{0}^{1,p(x)}(\Omega)$ which we call generalized Lebesgue-Sobolev
spaces. Firstly we state some basic properties of spaces
$W_{0}^{1,p(x)}(\Omega)$ which will be used later (for details see
\cite{17,14,15,4} and the references therein).

We write
\begin{align*}
&C_{+}(\overline{\Omega})=\{h|h\in C(\overline{\Omega}),h(x)>1\text{ for any }x\in\overline{\Omega}\},\\
&h^{+}=\mathop\text{max}\limits_{\overline{\Omega}}h(x),h^{-}=\mathop\text{min}\limits_{\overline{\Omega}}h(x)\text{ for any }h\in C({\overline{\Omega}}),\\
&L^{p(x)}(\Omega)=\{u|u\text{ is a measurable real-valued function, }\int_{\Omega}|u(x)|^{p(x)}\text{d}x<\infty\}.
\end{align*}
The linear vector space $L^{p(x)}(\Omega)$ can be equipped by the
following norm
\begin{align*}
|u|_{p(x)}=\text{inf}\bigg\{\lambda>0\bigg|\int_{\Omega}\bigg|\frac{u(x)}{\lambda}\bigg|^{p(x)}\text{d}x\leq1\bigg\},
\end{align*}
then $(L^{p(x)},|\cdot|_{p(x)})$ becomes a Banach space and we call it variable exponent Lebesgue space.\\

In the following we shall collect some basic propositions concerning the variable exponent Lebesgue spaces. These propositions will be used throughout our analysis.\\\\
\textbf{Proposition 1.1} (see Fan and Zhao \cite{17} and Zhao et al. \cite{14}).\emph{ (i) The space $(L^{p(x)},|\cdot|_{p(x)})$ is a separable, uniform convex Banach space, and its conjugate space is $L^{q(x)}(\Omega)$, where $\frac{1}{p(x)}+\frac{1}{q(x)}=1$. For any $u\in L^{p(x)}$ and $v\in L^{q(x)}$, we have
\begin{align*}
\bigg|\int_{\Omega}uv\text{\emph{d}}x\bigg|\leq(\frac{1}{p^{-}}+\frac{1}{q^{-}})|u|_{p(x)}|v|_{q(x)}.
\end{align*}
(ii) If $p_{1}, p_{2}\in C_{+}(\Omega), p_{1}(x)\leq p_{2}(x)$ for any $x\in\overline{\Omega}$, then $L^{p_{2}(x)}(\Omega)\hookrightarrow L^{p_{1}(x)}(\Omega)$, and the imbedding is continuous.
}\\\\
\textbf{Proposition 1.2} (see Fan and Zhao \cite{17} and Zhao et al. \cite{15}).\emph{ If $f:\Omega\times\mathbb{R}\rightarrow\mathbb{R}$ is a Caratheodory function and satisfies
\begin{align*}
|f(x,s)|\leq a(x)+b|s|^{\frac{p_{1}(x)}{p_{2}(x)}}\text{ for any }x\in\Omega,s\in\mathbb{R},
\end{align*}
where $p_{1},p_{2}\in C_{+}(\overline{\Omega}), a(x)\in L^{p_{2}(x)}(\Omega), a(x)\geq0$ and $b\geq0$ is a constant, then the Nemytsky operator from $L^{p_{1}(x)}(\Omega)$ to $L^{p_{2}(x)}(\Omega)$ defined by $(N_{f}(u))(x)=f(x,u(x))$ is a continuous and bounded operator.
}\\\\
\textbf{Proposition 1.3} (see Fan and Zhao \cite{17} and Zhao et al. \cite{14}).\emph{ If we denote
\begin{align*}
\rho(u)=\int_{\Omega}|u|^{p(x)}\text{\emph{d}}x, \forall u\in L^{p(x)}(\Omega),
\end{align*}
then\\
(i) $|u|_{p(x)}<1(=1;>1)\Leftrightarrow\rho(u)<1(=1;>1)$;\\
(ii) $|u|_{p(x)}>1\Rightarrow |u|_{p(x)}^{p^{-}}\leq\rho(u)\leq|u|_{p(x)}^{p^{+}};\,\,|u|_{p(x)}<1\Rightarrow |u|_{p(x)}^{p^{-}}\geq\rho(u)\geq|u|_{p(x)}^{p^{+}}$;\\
(iii)
$|u|_{p(x)}\rightarrow0\Leftrightarrow\rho(u)\rightarrow0;\,\,|u|_{p(x)}\rightarrow\infty\Leftrightarrow\rho(u)\rightarrow\infty.$
}\\\\
\textbf{Proposition 1.4} (see Fan and Zhao \cite{17} and Zhao et al. \cite{14}).\emph{ If $u, u_{n}\in L^{p(x)}(\Omega), n=1,2,...,$ then the following statements are equivalent to each other:\\
(1) $\text{\emph{lim}}_{k\rightarrow\infty}|u_{k}-u|_{p(x)}=0$;\\
(2) $\text{\emph{lim}}_{k\rightarrow\infty}\rho(u_{k}-u)=0$;\\
(3) $u_{k}\rightarrow u$ in measure in $\Omega$ and $\text{\emph{lim}}_{k\rightarrow\infty}\rho(u_{k})=\rho(u)$.
}\\\\
The variable exponent Sobolev space $W^{1,p(x)}(\Omega)$ is defined
by
\begin{align*}
W^{1,p(x)}(\Omega)=\{u\in L^{p(x)}(\Omega)| | \nabla u|\in L^{p(x)}(\Omega)\}
\end{align*}
and it is equipped with the norm
\begin{align*}
\|u\|_{p(x)}=|u|_{p(x)}+|\nabla u|_{p(x)},\forall u\in W^{1,p(x)}(\Omega).
\end{align*}
If we denote by $W_{0}^{1,p(x)}(\Omega)$ the closure of $C_{0}^{\infty}(\Omega)$ in $W^{1,p(x)}(\Omega)$ and
\begin{align*}
p^{*}(x)=\left\{ \begin{array}{rcl}
         &\frac{Np(x)}{N-p(x)}, &  p(x)<N, \\
         &\infty,&p(x)\geq N,
          \end{array}\right.
\end{align*}
we have the following\\\\
\textbf{Proposition 1.5} (see Fan and Zhao \cite{17}).\emph{
(i) $W^{1,p(x)}(\Omega)$ and $W_{0}^{1,p(x)}(\Omega)$ are separable reflexive Banach spaces;\\
(ii) If $q\in C_{+}(\overline{\Omega})$ and $q(x)<p^{*}(x)$ for any $x\in \overline{\Omega}$, then the imbedding from $W^{1,p(x)}(\Omega)$ to $L^{q(x)}(\Omega)$ is compact and continuous;\\
(iii) There is a constant $C>0$ such that
\begin{align*}
|u|_{p(x)}\leq C|\nabla u|_{p(x)}, \forall u\in W_{0}^{1,p(x)}(\Omega).
\end{align*}
}\\\\
\textbf{Remark 1.6.} \emph{
By (iii) of Proposition 1.5, we know that $|\nabla u|_{p(x)}$ and $\|u\|_{p(x)}$ are equivalent norms on $W_{0}^{1,p(x)}(\Omega)$. We will use $|\nabla u|_{p(x)}$ instead of $\|u\|_{p(x)}$ in the following discussions.
}\\

\section*{2. Properties of $(p_{1}(x),p_{2}(x))$-Laplace operator}

In this section we give the properties of the $(p_{1}(x),p_{2}(x))$-Laplace operator $(-\Delta_{p_{1}(x)}-\Delta_{p_{2}(x)})u:=-\text{div}(|\nabla u|^{p_{1}(x)-2}\nabla u)-\text{div}(|\nabla u|^{p_{2}(x)-2}\nabla u)$. Consider the following functional,
\begin{align*}
J(u)=\int_{\Omega}\frac{1}{p_{1}(x)}|\nabla u|^{p_{1}(x)}\text{d}x+\int_{\Omega}\frac{1}{p_{2}(x)}|\nabla u|^{p_{2}(x)}\text{d}x, \forall u\in X,
\end{align*}
where $X:=W_{0}^{1,p_{1}(x)}(\Omega)\cap W_{0}^{1,p_{2}(x)}(\Omega)$ with its norm given by $\|u\|:=\|u\|_{p_{1}(x)}+\|u\|_{p_{2}(x)},\forall u\in X$.

It is obvious that $J\in C^{1}(X,\mathbb{R})$ (see \cite{18}), and the $(p_{1}(x),p_{2}(x))$-Laplace operator is the derivative
operator of $J$ in the weak sense. Denote $L=J':X\rightarrow X^{*}$,
then
\begin{align*}
\langle L(u),v\rangle=\int_{\Omega}|\nabla u|^{p_{1}(x)-2}\nabla u\nabla v\text{d}x+\int_{\Omega}|\nabla u|^{p_{2}(x)-2}\nabla u\nabla v\text{d}x,\forall u,v\in X,
\end{align*}
in which $\langle\cdot,\cdot\rangle$ is the dual pair between $X$ and its dual $X^{*}$.\\\\
\textbf{Remark 2.1.} \emph{
$(X,\| \cdot \|)$ is a separable reflexive Banach space.
}\\\\\textbf{Theorem 2.2.} \emph{
(i) $L:X\rightarrow X^{*}$ is a continuous, bounded and strictly monotone operator;\\
(ii) $L$ is a mapping of type $(S_{+})$, namely: if $u_{n}\rightharpoonup u$ in $X$ and $\overline{\text{lim}}_{n\rightarrow\infty}\langle L(u_{n})-L(u),u_{n}-u\rangle\leq 0$, then $u_{n}\rightarrow u$ in $X$;\\
(iii) $L:X\rightarrow X^{*}$ is a homeomorphism.
}\\\\
\textbf{Proof.} (i) It is obvious that $L$ is continuous and
bounded. Considering the duality pair $\langle Lu-Lv,u-v\rangle$, we
have
\begin{align*}
       &\langle Lu-Lv,u-v\rangle\\
     =&\int_{\Omega}(|\nabla u|^{p_{1}(x)-2}\nabla u-|\nabla v|^{p_{1}(x)-2}\nabla v)(\nabla u-\nabla v)\text{d}x\\
       &+\int_{\Omega}(|\nabla u|^{p_{2}(x)-2}\nabla u-|\nabla v|^{p_{2}(x)-2}\nabla v)(\nabla u-\nabla v)\text{d}x\\
\geq&\int_{\Omega}h_{1}(x)\text{d}x+\int_{\Omega}h_{2}(x)\text{d}x\geq 0,
\end{align*}
where
\begin{align*}
0\leq h_{i}(x)=\left\{ \begin{array}{rcl}
         &\frac{(p_{i}(x)-1)|\nabla u- \nabla v|^{2}}{(|\nabla u|^{p_{i}(x)}+|\nabla v|^{p_{i}(x)})^{\frac{2-p_{i}(x)}{p_{i}(x)}}}, &  1<p_{i}(x)<2, \\
         &\frac{1}{2^{p_{i}(x)}}|\nabla u- \nabla v|^{p_{i}(x)},&p_{i}(x)\geq 2,
          \end{array}\right.
\end{align*}
$i=1,2$. Here we have applied the following inequalities (see
\cite{19,20}): for any $\xi,\eta\in \mathbb{R}^{N},$
\begin{align*}
&[(|\xi|^{p-2}\xi-|\eta|^{p-2}\eta)\cdot(\xi-\eta)](|\xi|^{p}+|\eta|^{p})^{\frac{2-p}{p}}\geq(p-1)|\xi-\eta|^{2}, &1<p<2;\\
&(|\xi|^{p-2}\xi-|\eta|^{p-2}\eta)\cdot(\xi-\eta)\geq\frac{1}{2^{p}}|\xi-\eta|^{p}, &p\geq 2,
\end{align*}
in which the equal-sign holds if and only if $\xi=\eta$.\\
(ii) From (i), if $u_{n}\rightharpoonup u$ in $X$ and $\overline{\text{lim}}_{n\rightarrow\infty}\langle L(u_{n})-L(u),u_{n}-u\rangle\leq 0$, then $\text{lim}_{n\rightarrow\infty}\langle L(u_{n})-L(u),u_{n}-u\rangle=0$. This implies by the proof of (i)
\begin{align*}
o_n(1)=\int_{\{1<p_{i}(x)<2\}}\frac{(p_{i}(x)-1)|\nabla u_{n}- \nabla u|^{2}}{(|\nabla u_{n}|^{p_{i}(x)}+|\nabla u|^{p_{i}(x)})^{\frac{2-p_{i}(x)}{p_{i}(x)}}}\text{d}x+\int_{\{p_{i}(x)\geq 2\}}\frac{1}{2^{p_{i}(x)}}|\nabla u_{n}- \nabla u|^{p_{i}(x)}\text{d}x.
\end{align*}
The above equality implies $\nabla u_{n}$ converges in measure to
$\nabla u$ in $\Omega$. So we may assume up to a subsequence (still
denoted by $\{\nabla u_{n}\}$) that $\nabla u_{n}(x)\rightarrow
\nabla u(x)$, a. e. $x\in\Omega$. By Fatou's Lemma, we get
\begin{align}
\mathop{\underline{\text{lim}}}\limits_{n\rightarrow\infty}\int_{\Omega}\frac{1}{p_{i}(x)}|\nabla
u_{n}|^{p_{i}(x)}
\text{d}x\geq\int_{\Omega}\frac{1}{p_{i}(x)}|\nabla
u|^{p_{i}(x)}\text{d}x.
\end{align}
In view of the assumption $u_{n}\rightharpoonup u$, we have
\begin{align}
\mathop{\text{lim}}\limits_{n\rightarrow\infty}\langle L(u_{n}),u_{n}-u\rangle=\mathop{\text{lim}}\limits_{n\rightarrow\infty}\langle L(u_{n})-L(u),u_{n}-u\rangle.
\end{align}
We also have
\begin{align*}
&\langle L(u_{n}),u_{n}-u\rangle\\
&=\int_{\Omega}|\nabla u_{n}|^{p_{1}(x)}\text{d}x-\int_{\Omega}|\nabla u_{n}|^{p_{1}(x)-2}\nabla u_{n}\nabla u\text{d}x+\int_{\Omega}|\nabla u_{n}|^{p_{2}(x)}\text{d}x-\int_{\Omega}|\nabla u_{n}|^{p_{2}(x)-2}\nabla u_{n}\nabla u\text{d}x\\
&\geq\int_{\Omega}|\nabla u_{n}|^{p_{1}(x)}\text{d}x-\int_{\Omega}|\nabla u_{n}|^{p_{1}(x)-1}|\nabla u|\text{d}x+\int_{\Omega}|\nabla u_{n}|^{p_{2}(x)}\text{d}x-\int_{\Omega}|\nabla u_{n}|^{p_{2}(x)-1}|\nabla u|\text{d}x\\
&\geq\int_{\Omega}|\nabla u_{n}|^{p_{1}(x)}\text{d}x-\int_{\Omega}(\frac{p_{1}(x)-1}{p_{1}(x)}|\nabla u_{n}|^{p_{1}(x)}+\frac{1}{p_{1}(x)}|\nabla u|^{p_{1}(x)})\text{d}x\\
&+\int_{\Omega}|\nabla u_{n}|^{p_{2}(x)}\text{d}x-\int_{\Omega}(\frac{p_{2}(x)-1}{p_{2}(x)}|\nabla u_{n}|^{p_{2}(x)}+\frac{1}{p_{2}(x)}|\nabla u|^{p_{2}(x)})\text{d}x\\
&\geq\int_{\Omega}\frac{1}{p_{1}(x)}|\nabla u_{n}|^{p_{1}(x)}\text{d}x-\int_{\Omega}\frac{1}{p_{1}(x)}|\nabla u|^{p_{1}(x)}\text{d}x+\int_{\Omega}\frac{1}{p_{2}(x)}|\nabla u_{n}|^{p_{2}(x)}\text{d}x-\int_{\Omega}\frac{1}{p_{2}(x)}|\nabla u|^{p_{2}(x)}\text{d}x.
\end{align*}
By (1) and (2), the above inequality implies
\begin{align}
\mathop{\text{lim}}\limits_{n\rightarrow\infty}\int_{\Omega}\frac{1}{p_{i}(x)}|\nabla u_{n}|^{p_{i}(x)}\text{d}x=\int_{\Omega}\frac{1}{p_{i}(x)}|\nabla u|^{p_{i}(x)}\text{d}x.
\end{align}
From (3) it follows that the integrals of the function family
$\{\frac{1}{p_{i}(x)}|\nabla u_{n}|^{p_{i}(x)}\}$ possesses
absolutely equicontinuity on $\Omega$ (see \cite{21}). Since
\begin{align*}
\frac{1}{p_{i}(x)}|\nabla u_{n}-\nabla u|^{p_{i}(x)}\leq C(\frac{1}{p_{i}(x)}|\nabla u_{n}|^{p_{i}(x)}+\frac{1}{p_{i}(x)}|\nabla u|^{p_{i}(x)}),
\end{align*}
the integrals of the family $\{\frac{1}{p_{i}(x)}|\nabla u_{n}-\nabla u|^{p_{i}(x)}\}$ are also absolutely equicontinuous on $\Omega$ and therefore
\begin{align}
\mathop{\text{lim}}\limits_{n\rightarrow\infty}\int_{\Omega}\frac{1}{p_{i}(x)}|\nabla u_{n}-\nabla u|^{p_{i}(x)}\text{d}x=0.
\end{align}
By (4),
\begin{align}
\mathop{\text{lim}}\limits_{n\rightarrow\infty}\int_{\Omega}|\nabla u_{n}-\nabla u|^{p_{i}(x)}\text{d}x=0.
\end{align}
From (5), we have $u_{n}\rightarrow u$ in $X$. This implies that $L$ is of type $(S_{+})$.\\
(iii) By (i), $L$ is an injection. Since
\begin{align*}
\mathop{\text{lim}}\limits_{\|u\|\rightarrow\infty}\frac{\langle Lu,u\rangle}{\|u\|}=\mathop{\text{lim}}\limits_{\|u\|\rightarrow\infty}\frac{\int_{\Omega}|\nabla u|^{p_{1}(x)}\text{d}x+\int_{\Omega}|\nabla u|^{p_{2}(x)}\text{d}x}{\|u\|}=\infty,
\end{align*}
$L$ is coercive, thus $L$ is a surjection by Minty-Browder Theorem (see \cite{22}). Hence $L$ has an inverse mapping $L^{-1}:X^{*}\rightarrow X$. In order to show that
 $L$ is a homeomorphism, it remains to show the continuity of $L^{-1}$.

For any $f_{n}, f\in X^{*}$ such that $f_{n}\rightarrow f$ in $X^{*}$, suppose $u_{n}=L^{-1}(f_{n}),u=L^{-1}(f)$, then we have $L(u_{n})=f_{n},L(u)=f$. Since $L$ is coercive, $\{u_{n}\}$ is bounded in X. We can assume that $u_{n_{k}}\rightharpoonup u_{0}$ in $X$. By $f_{n_{k}}\rightarrow f$ in $X^{*}$, we have
\begin{align*}
\mathop{\text{lim}}\limits_{n\rightarrow\infty}\langle
L(u_{n_{k}})-L(u_{0}),u_{n_{k}}-u_{0}\rangle=\mathop{\text{lim}}\limits_{n\rightarrow\infty}\langle
f_{n_{k}},u_{n_{k}}-u_{0}\rangle=\mathop{\text{lim}}\limits_{n\rightarrow\infty}\langle
f_{n_{k}}-f,u_{n_{k}}-u_{0}\rangle=0.
\end{align*}
Since $L$ is of type $S_{+}$, $u_{n_{k}}\rightarrow u_{0}$. By
injectivity of $L$, we have $u_{0}=u$. So $u_{n_{k}}\rightarrow u$.
We claim that $u_{n}\rightarrow u$ in $X$. Otherwise, there would
exist a subsequence $\{u_{m_{j}}\}$ of $\{u_{n}\}$ in $X$ and an
$\epsilon_{0}>0$, such that for any $m_{j}>0$, we have
$\|u_{m_{j}}-u\|\geq\epsilon_{0}$. But reasoning as above,
$\{u_{m_{j}}\}$ would contain a further subsequence
$u_{m_{j_{l}}}\rightarrow u$ in $X$ as $l\rightarrow\infty$, which
is a contradiction to $\|u_{m_{j_{l}}}-u\|\geq\epsilon_{0}$.
$\square$

\section*{3. Solutions to the equation}
In this section we will give the existence results of weak solutions
to problem $(P)$. We denote
\begin{align*}
p_{M}(x)=\text{max}\{p_{1}(x),p_{2}(x)\},\,p_{m}(x)=\text{min}\{p_{1}(x),p_{2}(x)\},\forall
x\in \Omega.
\end{align*}
It is easy to see that $p_{M}(x),p_{m}(x)\in
C_{+}(\overline{\Omega})$. For $q(x)\in C_+(\overline\Omega)$ such
that $q(x)<p_{M}^{*}(x)$ for any $x\in\overline\Omega$, we have
$X:=W_{0}^{1,p_{1}(x)}(\Omega)\cap
W_{0}^{1,p_{2}(x)}(\Omega)=W_{0}^{1,p_{M}(x)}(\Omega)\hookrightarrow
L^{q(x)}(\Omega)$, and the imbedding is continuous and compact. To
proceed, we give the definition of weak solution to problem $(P)$:\\

\noindent\textbf{Definition 3.1.}\emph{ We say that $u\in X$ is a
weak solution of $(P)$ if the following equality
\begin{align*}
\int_{\Omega}|\nabla u|^{p_{1}(x)-2}\nabla u\nabla v\text{d}x+\int_{\Omega}|\nabla u|^{p_{2}(x)-2}\nabla u\nabla v\text{d}x=\int_{\Omega}f(x,u)v\text{d}x
\end{align*}
holds for for any $v\in X:=W_{0}^{1,p_{1}(x)}(\Omega)\cap
W_{0}^{1,p_{2}(x)}(\Omega)$.
}\\\\
If $f$ is independent of $u$, we have\\\\
\textbf{Theorem 3.2.}\emph{
If $f(x,u)=f(x)$, and $f\in L^{\alpha(x)}(\Omega)$, where $\alpha\in C_{+}(\overline{\Omega})$, satisfies $\frac{1}{\alpha(x)}+\frac{1}{p^{*}_{M}(x)}<1$, then $(P)$ has a unique weak solution.
}\\\\
\textbf{Proof.}  Because $W_{0}^{1,p_{M}(x)}(\Omega)=W_{0}^{1,p_{1}(x)}(\Omega)\cap W_{0}^{1,p_{2}(x)}(\Omega)$, by Proposition 1.5 (ii), $\langle f, v\rangle:=\int_{\Omega}f(x)v\text{d}x$ (for any $v\in X$) defines a continuous linear functional on $X$. By Theorem 3.2 $L$ is a homeomorphism, $(P)$ has a unique solution. $\square$\\

From now on we always suppose that $f(x,t)$ satisfies assumption $(\text{f}_{0})$:\\
$(\text{f}_{0})$ $f:\Omega\times\mathbb{R}\rightarrow\mathbb{R}$ satisfies Caratheodory condition and
\begin{align*}
|f(x,t)|\leq C_{1}+C_{2}|t|^{\alpha(x)-1},\forall(x,t)\in\Omega\times\mathbb{R},
\end{align*}
where $\alpha\in C{+}(\overline{\Omega})$ and
$\alpha(x)<p^{*}_{M}(x)$ for any $x\in\overline\Omega$.

Let $\varphi(u):=\int_{\Omega}\frac{1}{p_{1}(x)}|\nabla
u|^{p_{1}(x)}\text{d}x+\int_{\Omega}\frac{1}{p_{2}(x)}|\nabla
u|^{p_{2}(x)}\text{d}x-\int_{\Omega}F(x,u)\text{d}x, u\in X$, where
$F(x,t)=\int^{t}_{0}f(x,s)ds$. It is trivial that $\varphi\in
C^{1}(X,\mathbb{R})$. Therefore, weak solutions of $(P)$ correspond to critical points of $\varphi$.

Let $g(u):=\int_{\Omega}F(x,u)\text{d}x$. Then $g':X\rightarrow X^{*}$ is completely continuous,
i.e. $u_{n}\rightharpoonup u$ in $X$ implies $g'(u_{n})\rightarrow g'(u)$ in $X^{*}$, thus $g$ is weakly continuous.\\\\
\textbf{Theorem 3.3.}\emph{ If $f$ satisfies the following
condition,
\begin{align}
|f(x,t)|\leq C_{1}+C_{2}|t|^{\beta(x)-1},\text{ where } 1\leq\beta^{+}<p_{m}^{-},
\end{align}
then $(P)$ has a weak solution.
}\\\\
\textbf{Proof.} By Condition (6) we have $|F(x,t)|\leq C(1+|t|^{\beta(x)})$, and
\begin{align*}
\varphi(u)&=\int_{\Omega}\frac{1}{p_{1}(x)}|\nabla u|^{p_{1}(x)}\text{d}x+\int_{\Omega}\frac{1}{p_{2}(x)}|\nabla u|^{p_{2}(x)}\text{d}x-\int_{\Omega}F(x,u)\text{d}x\\
&\geq \frac{1}{p_{M}^{+}}(\|u_{n}\|_{p_{1}(x)}^{p_{1}^{-}}+\|u_{n}\|_{p_{2}(x)}^{p_{2}^{-}})-C\int_{\Omega}|u|^{\beta(x)}\text{d}x-C_{3}\\
&\geq \frac{C_{5}}{p_{M}^{+}}\|u\|^{p_{m}^{-}}-C_{4}\|u\|^{\beta^{+}}-C_{3}\rightarrow\infty,\text{ as }\|u\|\rightarrow\infty.
\end{align*}
In view that $\varphi$ is also weakly lower semicontinuous, we see $\varphi$ has a global minimum point $u\in X$, which is a weak solution to problem $(P)$. We now complete the proof. $\square$\\\\
\textbf{Definition 3.4.} \emph{We say that the function $\varphi\in C^{1}(X,\mathbb{R})$ satisfies the Palais-Smale $(PS)$ condition in $X$ if any sequence $\{u_{n}\}\in X$ such that
\begin{align*}
|\varphi(u_{n})|\leq B,\text{ for some }B\in\mathbb{R};\\
\varphi'(u_{n})\rightarrow 0\text{ in
}X^{*}\text{ as }n\rightarrow\infty
\end{align*}
has a convergent subsequence.}\\\\
\textbf{Lemma 3.5.} \emph{
If $f$ satisfies\\
$(\text{\emph{f}}_{1})$ $\exists M>0, \theta>p_{M}^{+}$ such that
\begin{align*}
0<\theta F(x,t)\leq tf(x,t), |t|\geq M, x\in\Omega,
\end{align*}
then $\varphi$ satisfies (PS) condition.}\\\\
\textbf{Proof.} Suppose that $\{u_{n}\}\subset X$, $\{\varphi(u_{n})\}$ is bounded and $\|\varphi'(u_{n})\|\rightarrow0$. Without loss of generality, we could assume $\|u_n \|\geq M\geq 1$ for some large constant $M$, and then we have
\begin{align*}
C&\geq\varphi(u_{n})=\int_{\Omega}\frac{1}{p_{1}(x)}|\nabla u_{n}|^{p_{1}(x)}\text{d}x+\int_{\Omega}\frac{1}{p_{2}(x)}|\nabla u_{n}|^{p_{2}(x)}\text{d}x-\int_{\Omega}F(x,u_{n})\text{d}x\\
&\geq \int_{\Omega}\frac{1}{p_{1}(x)}|\nabla u_{n}|^{p_{1}(x)}\text{d}x+\int_{\Omega}\frac{1}{p_{2}(x)}|\nabla u_{n}|^{p_{2}(x)}\text{d}x-\int_{\Omega}\frac{u_{n}}{\theta}f(x,u_{n})\text{d}x-c\\
&\geq \int_{\Omega}(\frac{1}{p_{1}(x)}-\frac{1}{\theta})|\nabla u_{n}|^{p_{1}(x)}\text{d}x+\int_{\Omega}(\frac{1}{p_{2}(x)}-\frac{1}{\theta})|\nabla u_{n}|^{p_{2}(x)}\text{d}x\\
&+\int_{\Omega}\frac{1}{\theta}(|\nabla u_{n}|^{p_{1}(x)}+|\nabla u_{n}|^{p_{2}(x)}-u_{n}f(x,u_{n}))\text{d}x-c\\
&\geq (\frac{1}{p_{M}^{+}}-\frac{1}{\theta})(\|u_{n}\|_{p_{1}(x)}^{p_{1}^{-}}+\|u_{n}\|_{p_{2}(x)}^{p_{2}^{-}})+\int_{\Omega}\frac{1}{\theta}(|\nabla u_{n}|^{p_{1}(x)}+|\nabla u_{n}|^{p_{2}(x)}-u_{n}f(x,u_{n}))\text{d}x-c\\
&\geq C_{5}(\frac{1}{p_{M}^{+}}-\frac{1}{\theta})\|u_{n}\|^{p_{m}^{-}}-\frac{1}{\theta}\|\varphi'(u_{n})\|\|u_{n}\|-c.
\end{align*}
The above inequality implies $\{u_{n}\}$ is bounded in $X$. Without loss of generality, assume that $u_{n}\rightharpoonup u$, then $g'(u_{n})\rightarrow g'(u)$. Since $\varphi'(u_{n})=L(u_{n})-g'(u_{n})$, we can get $L(u_{n})\rightarrow g'(u)$. Since $L$ is a homeomorphism, we conclude $u_{n}\rightarrow u$, which implies $\varphi$ satisfies (PS) condition. $\square$\\\\
\textbf{Theorem 3.6.} \emph{
If $f$ satisfies $(\text{\emph{f}}_{0})$, in which $\alpha^{-}>p_{M}^{+}$, $(\text{\emph{f}}_{1})$ and the following\\
$(\text{\emph{f}}_{2})$ $f(x,t)=o(|t|^{p_{M}^{+}-1}),t\rightarrow 0$ for $x\in \Omega$ uniformly,\\
then $(P)$ has a nontrivial solution.}\\\\
\textbf{Proof.} We shall show that $\varphi$ satisfies conditions of
Mountain Pass theorem. \\(i) From Lemma 3.5, $\varphi$ satisfies
$(PS)$ condition in $X$. Since
$p^{+}<\alpha^{-}\leq\alpha(x)<p^{*}(x)$, $X\hookrightarrow
L^{p_{M}^{+}}(\Omega)$, i.e. there exists $C_{0}>0$ such that
\begin{align*}
|u|_{p^{+}}\leq C_{0}\|u\|,\forall u\in X.
\end{align*}
Let $\epsilon>0$ be small enough such that $\epsilon C_{0}^{p_{M}^{+}}\leq\frac{C_{6}}{2p_{M}^{+}}$. By the assumptions $\text{f}_{0}$ and $\text{f}_{2}$, we have
\begin{align*}
F(x,t)\leq \epsilon|t|^{p_{M}^{+}}+C(\epsilon)|t|^{\alpha(x)}\forall(x,t)\in\Omega\times\mathbb{R}.
\end{align*}
For $\|u\|\leq 1$ we have the following
\begin{align*}
\varphi(u)&=\int_{\Omega}\frac{1}{p_{1}(x)}|\nabla u|^{p_{1}(x)}\text{d}x+\int_{\Omega}\frac{1}{p_{2}(x)}|\nabla u|^{p_{2}(x)}\text{d}x-\int_{\Omega}F(x,u)\text{d}x\\
&\geq \frac{1}{p_{M}^{+}}(\|u\|_{p_{1}(x)}^{p_{1}^{+}}+\|u\|_{p_{2}(x)}^{p_{2}^{+}})-\epsilon\int_{\Omega}|u|^{p_{M}^{+}}\text{d}x-C(\epsilon)\int_{\Omega}|u|^{\alpha(x)}\text{d}x\\
&\geq \frac{C_{6}}{p_{M}^{+}}\|u\|^{p_{M}^{+}}-\epsilon C_{0}^{p_{M}^{+}}\|u\|^{p_{M}^{+}}-C(\epsilon)\|u\|^{\alpha^{-}}\\
&\geq \frac{C_{6}}{2p_{M}^{+}}\|u\|^{p_{M}^{+}}-C(\epsilon)\|u\|^{\alpha^{-}},
\end{align*}
which implies the existence of $r\in (0,\,1)$ and $\delta>0$ such that $\varphi(u)\geq\delta>0$ for every $u\in X$ satisfies $\|u\|=r$.\\
(ii) From $(\text{f}_{1})$ we see for suitable positive constants $C, C'$,
\begin{align*}
F(x,t)\geq C|t|^{\theta}-C',\forall x\in\overline{\Omega}, t\geq M.
\end{align*}
For any fixed $w\in X\backslash\{0\}$, and $t>1$, we have
\begin{align*}
\varphi(tw)&=\int_{\Omega}\frac{1}{p_{1}(x)}|t\nabla w|^{p_{1}(x)}\text{d}x+\int_{\Omega}\frac{1}{p_{2}(x)}|t\nabla w|^{p_{2}(x)}\text{d}x-\int_{\Omega}F(x,tw)\text{d}x\\
&\leq t^{p_{M}^{+}}(\int_{\Omega}\frac{1}{p_{1}(x)}|\nabla w|^{p_{1}(x)}\text{d}x+\int_{\Omega}\frac{1}{p_{2}(x)}|\nabla w|^{p_{2}(x)}\text{d}x)-Ct^{\theta}\int_{\Omega}|w|^{\theta}\text{d}x-C_{7}\\
&\rightarrow-\infty, \text{ as } t\rightarrow+\infty.
\end{align*}
(iii) It is obvious $\varphi(0)=0$.

From (i), (ii) and (iii), we conclude $\varphi$ satisfies the conditions of Mountain Pass Theorem (see \cite{18}). So $\varphi$ admits at least one nontrivial critical point. $\square$\\\\
\textbf{Theorem 3.7.} \emph{
If $f$ satisfies $(\text{\emph{f}}_{0})$, $(\text{\emph{f}}_{1})$, $p_{M}^{+}<\alpha^{+}$ and the following\\
$(\text{\emph{f}}_{3})$ $f(x,-t)=-f(x,t)$, $x\in \Omega,t\in\mathbb{R}$,\\
then $\varphi$ has a sequence of critical points $\{u_{n}\}$ such that $\varphi(u_{n})\rightarrow+\infty$ and $(P)$ has infinite many pairs of solutions.}\\\\
Because $X$ is a reflexive and separable Banach space, there are $\{e_{j}\}\subset X$ and $\{e^{*}_{j}\}\subset X^{*}$ such that
\begin{align*}
X=\overline{\text{span}\{e_{j}|j=1,2,...\}},X^{*}=\overline{\text{span}\{e^{*}_{j}|j=1,2,...\}}
\end{align*}
and
\begin{align*}
\langle e_{i}^{*}, e_{j}\rangle=\left\{ \begin{array}{rcl}
         &1,&  i=j; \\
         &0,&i\neq j.
          \end{array}\right.
\end{align*}
Denote $X_{j}=\text{span}\{e_{j}\},Y_{k}=\oplus_{j=1}^{k}X_{j},Z_{k}=\oplus_{j=k}^{\infty}X_{j}$.\\\\
\textbf{Lemma 3.8.} \emph{
If $\alpha\in C_{+}(\overline{\Omega}),\alpha(x)<p_{M}^{*}(x)$ for any $x\in\overline{\Omega}$, denote
\begin{align*}
\beta_{k}=\text{sup}\{|u|_{\alpha(x)}|\|u\|=1,u\in Z_{k}\},
\end{align*}
then $\text{\emph{lim}}_{k\rightarrow\infty}\beta_{k}=0$.}\\\\
\textbf{Proof.} For $0<\beta_{k+1}\leq\beta_{k}$, then $\beta_{k}\rightarrow\beta\geq0$. Suppose $u_{k}\in Z_{k}$ satisfy
\begin{align*}
\|u_{k}\|=1,0\leq\beta_{k}-|u_{k}|_{\alpha(x)}<\frac{1}{k},
\end{align*}
then we may assume up a subsequence that $u_{k}\rightharpoonup u$ in
$X$, and
\begin{align*}
\langle e^{*}_{j},u\rangle=\text{lim}_{k\rightarrow\infty}\langle e^{*}_{j},u_{k}\rangle=0,j=1,2,...,
\end{align*}
in which the last equal sign holds for $u_{k}\in Z_{k}$. The above equality implies that $u=0$, so $u_{k}\rightharpoonup 0$ in $X$. Since the imbedding from $X$ to $L^{\alpha(x)}(\Omega)$ is compact, we have $u_{k}\rightarrow 0$ in $L^{\alpha(x)}(\Omega)$. We finally get $\beta_{k}\rightarrow0$. $\square$\\\\
\textbf{Proof of Theorem 3.7.} By $(\text{f}_{1}),(\text{f}_{3})$, $\varphi$ is an even functional and satisfies $(PS)$ condition. We only need to prove that if $k$ is large enough, then there exist $\rho_{k}>\gamma_{k}>0$ such that\\
$(A_{1})$ $b_{k}:=\text{inf}\{\varphi(u)|u\in Z_{k},\|u\|=\gamma_{k}\}\rightarrow\infty,(k\rightarrow\infty)$;\\
$(A_{2})$ $a_{k}:=\text{max}\{\varphi(u)|u\in Y_{k},\|u\|=\rho_{k}\}\leq0$,\\
Then the conclusion of Theorem 3.7 can be obtained by Fountain Theorem (see \cite{16}) and Lemma 3.5.\\
$(A_{1})$ For any $u\in Z_{k}$ with $\|u\|$ is big enough, we have
\begin{align*}
\varphi(u)&=\int_{\Omega}\frac{1}{p_{1}(x)}|\nabla u|^{p_{1}(x)}\text{d}x+\int_{\Omega}\frac{1}{p_{2}(x)}|\nabla u|^{p_{2}(x)}\text{d}x-\int_{\Omega}F(x,u)\text{d}x\\
&\geq \frac{1}{p_{M}^{+}}(\|u\|_{p_{1}(x)}^{p_{1}^{-}}+\|u\|_{p_{2}(x)}^{p_{2}^{-}})-c\int_{\Omega}|u|^{\alpha(x)}\text{d}x-c_{1}\\
&\geq \frac{C_{5}}{p_{M}^{+}}\|u\|^{p_{m}^{-}}-c|u|_{\alpha(x)}^{\alpha(\xi)}-c_{2},\text{ where }\xi\in\Omega\\
&\geq\left\{ \begin{array}{rcl}
         &\frac{C_{5}}{p_{M}^{+}}\|u\|^{p_{m}^{-}}-c-c_{2},&\text{ if }|u|_{\alpha(x)}\leq 1; \\
         &\frac{C_{5}}{p_{M}^{+}}\|u\|^{p_{m}^{-}}-c\beta_{k}^{\alpha^{+}}\|u\|^{\alpha^{+}}-c_{2},&\text{ if }|u|_{\alpha(x)}> 1,
          \end{array}\right.\\
&\geq \frac{C_{5}}{p_{M}^{+}}\|u\|^{p_{m}^{-}}-c\beta_{k}^{\alpha^{+}}\|u\|^{\alpha^{+}}-c_{3}\\
&=
C_{5}(\frac{1}{p_{M}^{+}}\|u\|^{p_{m}^{-}}-c'\beta_{k}^{\alpha^{+}}\|u\|^{\alpha^{+}})-c_{3}.
\end{align*}
Set
$\|u\|=\gamma_{k}=(c'\alpha^{+}\beta_{k}^{\alpha^{+}})^{\frac{1}{p_{m}^{-}-\alpha^{+}}}$.
Because $\beta_{k}\rightarrow 0$ and $p_{m}^{-}\leq
p_{M}^{+}<\alpha^{+}$, we have
\begin{align*}
\varphi(u)&\geq C_{5}(\frac{1}{p_{M}^{+}}\|u\|^{p_{m}^{-}}-c'\beta_{k}^{\alpha^{+}}\|u\|^{\alpha^{+}})-c_{3}\\
&=C_{5}(\frac{1}{p_{M}^{+}}(c'\alpha^{+}\beta_{k}^{\alpha^{+}})^{\frac{p_{m}^{-}}{p_{m}^{-}-\alpha^{+}}}-c'\beta_{k}^{\alpha^{+}}(c'\alpha^{+}\beta_{k}^{\alpha^{+}})^{\frac{\alpha^{+}}{p_{m}^{-}-\alpha^{+}}})-c_{3}\\
&=C_{5}(\frac{1}{p_{M}^{+}}-\frac{1}{\alpha^{+}})(c'\alpha^{+}\beta_{k}^{\alpha^{+}})^{\frac{p_{m}^{-}}{p_{m}^{-}-\alpha^{+}}}-c_{3}\rightarrow\infty, \text{ as }k\rightarrow\infty.
\end{align*}
$(A_{2})$ From $(\text{f}_{1})$, we get
\begin{align*}
F(x,t)\geq c_{1}|t|^{\theta}-c_{2}.
\end{align*}
Because $\theta>p_{M}^{+}$ and $\text{dim}Y_{k}=k$, it is easy to get $\varphi(u)\rightarrow-\infty$ as $\|u\|\rightarrow\infty$ for $u\in Y_{k}$. $\square$\\\\
\textbf{Remark 3.9.} \emph{ We can extend all our results to the $n$-case: Consider the following $(p_{1}(x), p_{2}(x),...,p_{n}(x))$-Laplace equation
\begin{align*}
(P_{n})\left\{ \begin{array}{rcl}
         &-\text{\emph{div}}(|\nabla u|^{p_{1}(x)-2}\nabla u)-\text{\emph{div}}(|\nabla u|^{p_{2}(x)-2}\nabla u)-...-\text{\emph{div}}(|\nabla u|^{p_{n}(x)-2}\nabla u)=f(x,u), &  \mbox{ in }
         \Omega; \\
         &u=0,&\mbox{ on } \partial\Omega.
          \end{array}\right.
\end{align*}
}\\\\

\renewcommand{\baselinestretch}{0.1}

\end{document}